\documentclass[11pt]{article}
\usepackage{amsmath}
\numberwithin{equation}{section}
\usepackage{amsfonts,amssymb,latexsym,amsbsy, bbm, theorem, mathrsfs, enumerate,color}
\usepackage{graphicx, graphics}
\usepackage{float,color,fancybox,shapepar,setspace,hyperref}
\usepackage{subfigure}
\usepackage{pgf,tikz}
\usetikzlibrary{arrows}
\newtheorem{Main Theorem}{Main Theorem}

\newtheorem{Conjecture}{Conjecture}
\newtheorem{Theorem}{Theorem}

\newtheorem{Lemma}{Lemma}
\newtheorem{Claim}{Claim}
\newtheorem{corollary}{Corollary}

\def\pf{\medskip\noindent {\emph{\bf Proof}.}~~}

\def\ora{\overrightarrow}
\def\ola{\overleftarrow}
\begin{document}

\title{Generalized Tur\'{a}n number  for linear forests}

\author{Xiutao ZHU and Yaojun CHEN\footnote{Corresponding author. Email:yaojunc@nju.edu.cn} \\
 \small{Department of Mathematics, Nanjing University, Nanjing 210093, P.R. CHINA}}
\date{}
\maketitle
\begin{abstract}

The generalized Tur\'{a}n number $ex(n,K_s,H)$ is defined to be the maximum number of copies of a complete graph $K_s$ in any $H$-free graph on $n$ vertices. Let $F$ be a linear forest consisting of $k$ paths of orders $\ell_1,\ell_2,...,\ell_k$. In this paper,  by characterizing the structure of the $F$-free graph with large minimum degree, we determine the value of $ex(n,K_s,F)$ for $n=\Omega\left(|F|^s\right)$ and $k\geq 2$ except some $\ell_i=3$, and the corresponding extremal graphs. The special case when $s=2$ of our result improves some results of Bushaw and Kettle (2011) and Lidick\'{y} et al. (2013) on the classical Tur\'{a}n number for linear forests. 

\vskip 2mm


\end{abstract}

\section{Introduction}
Let $G=(V(G),E(G))$ be a graph and $|G|=|V(G)|$ is the order of $G$. If $S\subseteq V(G)$, we use $G-S$ to denote the subgraph obtained from $G$ by deleting all vertices in $S$ and $G[S]$ the subgraph induced by $S$. Let $N_S(v)$ be the set of the neighbors of $v$ in $S$ and $d_S(v)=\left|N_S(v)\right|$, abbreviated as $N(v)$ and $d(v)$ for $S=V(G)$. We use $\delta(G)$, $\alpha(G)$, $\kappa(G)$ to denote the minimum degree, independence number, and connectivity of $G$, respectively.
Let $G\cup H$ denote the disjoint union of $G$ and $H$, and $kG$ denotes $k$ disjoint copies of $G$. We use $G+H$ to denote the graph obtained from $G\cup H$ by adding all edges between $V(G)$ and $V(H)$. A path, clique and empty graph on $m$ vertices are denoted by $P_m$, $K_m$ and $E_m$, respectively. If $u$ and $v$ are two vertices of a path $P$, we use $uPv$ to denote the segment of $P$ from $u$ to $v$. We denote by $\overrightarrow{P}$ the path $P$ with a given direction and $\overleftarrow{P}$ the path with the inverse direction.  We use $p(G)$ to denote the order of a longest path in $G$ and $P$ is hamiltonian if $|P|=|G|$. A path $P$ is called a strong dominating path if $N(v)\subseteq V(P)$ for any $v\in V(G)\backslash V(P)$.

A linear forest is the union of vertex disjoint paths, denoted by $F=P_{\ell_1}\cup \cdots \cup P_{\ell_k}$, where $\ell_1\ge\cdots\ge \ell_k$. Let $\mathcal{O}=\{\ell_i~|~\ell_i~is~odd,~1\leq i\leq k\}$ and $\chi_{\mathcal{O}}(\ell_i)$ denote the indicative function. Define 
$$\delta_F=\sum_{i=1}^k\left\lfloor\ell_i/2\right\rfloor-1$$
and $G_F(n)=K_{\delta_F}+E_{n-\delta_F}$ or $K_{\delta_F}+(E_{n-\delta_F-2}\cup K_2)$, according to $F$ contains an even component or not,  as shown in Figure 1($a$) and ($b$), respectively.

\begin{picture}(70,80)(10,8)
\centering
{\includegraphics[width=5.5in]{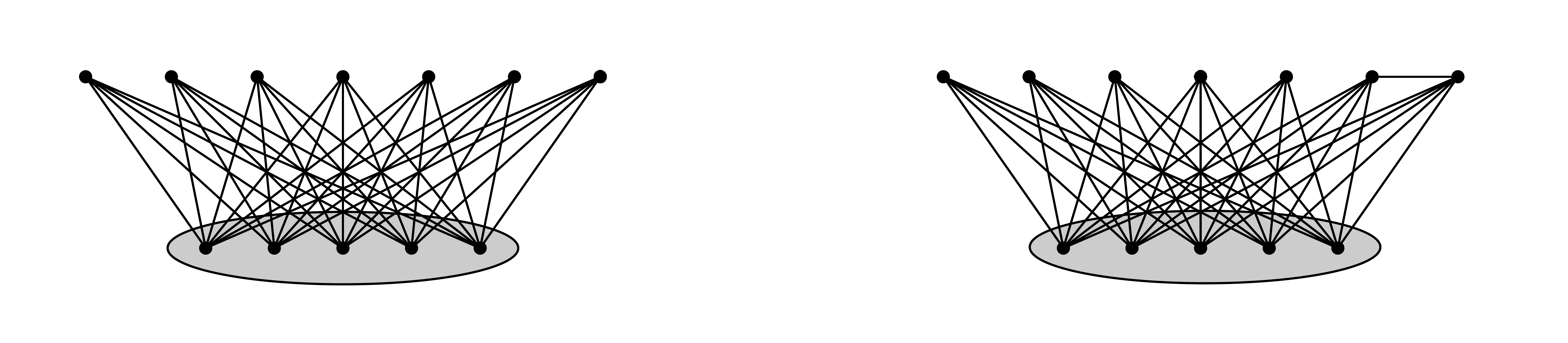}}
\put(-312,0){\makebox(3,3){($a$)}}
\put(-94,0){\makebox(3,3){($b$)}}
\put(-203,-20){\makebox(3,3){Figure 1. $G_F(n)$ (oval denotes a clique $K_{\delta_F}$)}}
\end{picture}
\vskip 12mm

The Tur\'an number of $H$, denoted by $ex(n,H)$, is the maximum number of edges in an $H$-free graph on $n$ vertices. Tur\'an theory dates back to the early 1940's when Tur\'an \cite{turan} determined the value of $ex(n, K_r)$ and extremal graph. For a given graph $H$ with chromatic number $\chi\geq 3$, an asymptotic solution to $ex(n,H)$ is given by the celebrated Erd\H{o}s-Stone-Simonovits theorem \cite{Er,stone} as below.
$$ex(n,H)=\left(\frac{\chi-2}{\chi-1}+o(1)\right){n\choose 2}.$$
In viewing this result, one can see that the asymptotic behavior of $ex(n,H)$ is still unknown for $H$ being a bipartite graph. The following is a classical result on $ex(n,H)$ when $H$ is bipartite due to Erd\H{o}s and Gallai.

\begin{Theorem}(Erd\H{o}s and Gallai \cite{Gallai})\label{Gallai} $ex(n,P_\ell)\le \frac{(\ell-2)n}{2}$, equality holds if and only if $(\ell-1)|n$.
\end{Theorem}
In the same paper, they also considered the graph $H=kK_2$ and determined $ex(n,kK_2)$ and the extremal graphs.
Note that both a path $P_\ell$ and $k$ independent edges $kK_2$ are special linear forests.  In 2011, Bushaw and Kettle \cite{Kettle} generalized the above results by determining the Tur\'an number and extremal graph for the linear forest consisting of the paths of the same lengths.
\begin{Theorem}(Bushaw and Kettle \cite{Kettle})\label{BK}
Let $F=kP_{\ell}$ be a linear forest so that $\ell\neq 3$ and $k\ge 2$. If $n= \Omega\left(|F|^{\frac{3}{2}}2^{\ell}\right)$, then
$$ex(n,F)=\binom{\delta_F}{2}+(n-\delta_F)\delta_F+\mu,$$
where $\mu=1$ if $\ell$ is odd and $\mu=0$ otherwise. The only extremal graph is $G_F(n)$.
\end{Theorem}
One can see that Theorem \ref{BK} is valid only when $n$ is exponential in $\ell$, which is very large. 
Furthermore, Bushaw and Kettle \cite{Kettle} posed the following.
\begin{Conjecture}(Bushaw and Kettle \cite{Kettle})\label{BKC}
Theorem \ref{BK} holds for  all $n=\Omega\left(|F|\right)$.
\end{Conjecture}
This conjecture is still open now. 
Soon after, Lidick\'y et al. \cite{Liu} extended Bushaw and Kettle's result by determining $ex(n,F)$ for $F\not= kP_3$ and sufficiently large $n$, and $G_F(n)$ is the extremal graph. Recently, Yuan and Zhang \cite{yuan} determined $ex(n,F)$  and the extremal graphs for all $n$ when $F$ contains at most one odd components.

In 2016, Alon and Shikhelman \cite{Alon} extended the idea of the Tur\'an number to the generalized Tur\'an number $ex(n,T,H)$, that is,  the maximum possible number of copies of $T$ in an $H$-free graph $G$ on $n$ vertices. Obviously, $ex(n,K_2,H)=ex(n,H)$. Before this concept was proposed formally, there are several papers dealing with the function $ex(n,T,H)$ for $T\neq K_2$, related papers can refer to \cite{Bollobas,Erdos,G,Gyori,H,Zykov}.
Since 2016, the problem of estimating generalized Tur\'{a}n number has received a lot of attentions. The most attractive problem maybe is to determine $ex(n,K_s, H)$.

Alon and Shikhelman \cite {Alon} generalized Erd\H{o}s-Stone-Simonovits theorem as below:
 $$ex(n,K_s,H)=\binom{\chi-1}{s}\left(\frac{n}{\chi-1}\right)^s+o(n^s),$$
where $\chi\geq 3$ is the chromatic number of $H$. Ma and Qiu \cite{Ma} improved the error term $o(n^s)$ up to a constant factor. Luo \cite{Luo} obtained the upper bound of $ex(n,K_s,P_\ell)$, which extends Theorem \ref{Gallai}.
\begin{Theorem} \label{Luo}(Luo \cite{Luo})
$ex(n,K_s,P_\ell)\le \frac{n}{\ell-1}\binom{\ell-1}{s}$ and the equality holds if and only if $(\ell-1)|n$.
\end{Theorem}
Wang \cite{Wang} determined $ex(n,K_s,kK_2)$, which extends the result obtained by Erd\H{o}s and Gallai in \cite{Gallai}.
Some other results about generalized Tur\'an number can be found in \cite{chase,A,k.G,peng}.

\vskip 5mm
Motivated by the results of Bushaw and Kettle \cite{Kettle} and Lidick\'y et al. \cite{Liu}, we consider the generalized Tur\'{a}n number $ex(n,K_s,H)$ for $H$ being a linear forest in more general situation. 
The main  result of this paper is the following.
\begin{Theorem}\label{Thm2}
Let $F$ be a linear forest so that  $\ell_i\not=3$ for $1\leq i\leq k$ and $k\ge 2$. Suppose $n\ge 5\binom{|F|-1}{s}^2/\left({(|F|-1)\binom{\delta_F}{s-1}}\right)+\delta_F$ and $s\le \delta_F+1$. 
If $F\neq 2P_\ell$ with $\ell$  odd, then 
 $$ex(n,K_s,F)=max\left\{\mathcal{N}_s(G_F(n)), ex(n,K_s,P_{\ell_1})\right\},$$ 
 and if $F=2P_\ell$ with $\ell$  odd, then 
\[
ex(n,K_s,F)=\left\{
\begin{array}{ll}
\mathcal{N}_s(G_F(n)),&~for~s\leq \delta_F, \\

\lfloor\frac{n-1}{\ell-1}\rfloor\ell+\mu, &~for~s=\delta_F+1,
\end{array}
\right.
\]
where $\mu=1$ if  $n-1\equiv \ell-2$ (mod $\ell-1$),  and $\mu=0$  otherwise.
\end{Theorem}

It is worth noting that in Theorem \ref{Thm2}, $\mathcal{N}_s(G_F(n))\ge ex(n,K_s,P_{\ell_1})$ is not always true. For example, if $s\ge3$ and $\ell_1\gg \ell_2+\cdots+\ell_k$,  then $\mathcal{N}_s(G_F(n))\le ex(n,K_s,P_{\ell_1})$. Thus it is reasonable to use $max\{\mathcal{N}_s(G_F(n)), ex(n,K_s,P_{\ell_1})\}$ rather than $\mathcal{N}_s(G_F(n))$.
 However, it is always true that $\mathcal{N}_s(G_F(n))\ge ex(n,K_s,P_{\ell_1})$ if $s=2$. Thus, we have
\begin{corollary}
Let $F$ be a linear forest so that  $\ell_i\not=3$ for $1\leq i\leq k$ and $k\ge 2$. If $n=\Omega\left(|F|^2\right)$, then
$ex(n,F)=e(G_F(n))$.
\end{corollary}
This implies Conjecture \ref{BKC} is indeed true for all $n=\Omega\left(|F|^2\right)$, which reduces the lower bound for $n$ in Theorem \ref{BK} from  exponential to polynomial in $|F|$.

\vskip 3mm
The main idea to prove Theorem \ref{Thm2} is to delete one vertex of small degree each step from an $F$-free graph until the remaining graph is empty or has large minimum degree. And then we count the number of copies of $K_s$ in the remaining graph. To do this, we have to characterize the structure of an $F$-free graph with high minimum degree, which is a key ingredient in proving  Theorem \ref{Thm2}.

\begin{Theorem}\label{Thm1}
Let $F$ be a linear forest so that  $\ell_i\not=3$ for $1\leq i\leq k$ and $k\ge2$. If $G$ is an $F$-free connected graph with $\delta(G)\ge \delta_F$ and $|G|=n\ge 2|F|$, then
$G\subseteq G_F(n)$ unless
\begin{enumerate}
 \item $F\in\{P_{\ell+1}\cup P_\ell$, $2P_{\ell}\}$ with $\ell$ even, or $F\in\{P_{\ell+2}\cup P_\ell$, $P_{\ell+1}\cup P_\ell, 2P_\ell\cup P_2\}$ with $\ell$ odd, and $G=K_1+tK_{\ell-1}$.
 \item $F=2P_\ell$ with $\ell$ odd, and $G\subseteq K_1+(t_1K_{\ell-1}\cup t_2K_{\ell-2})$.
 \item $F=3P_5$ and $G=K_2+tK_4$ or $G=E_2+tK_4$.

\end{enumerate}
\end{Theorem}

Theorem \ref{Thm1} is of interest in its own right. In fact,  Johansson \cite{Johansson} proved that: For any $F$, if $|G|=|F|$ and $\delta(G)\ge \delta_F+1$, then $F\subseteq G$. Later, Chen et al. \cite{chen} and  Egawa and Ota \cite{Egawa} independently obtained a stronger result: For any $F$, if $|G|\ge |F|$ and $\sigma_3(G)\ge 3(\delta_F+1)$, then $F\subseteq G$, where $\sigma_3(G)$ is the minimum degree sum of three nonadjacent vertices.   Chen and Zhang \cite{zhang} obtained a version of stability result for $F$ contains at most one odd path, they characterized the structure of connected graph $G$ if $G$ contains no $F$ and $\delta(G)=\delta_F$.

The remainder of this paper is organized as follows. In Section 2, we give some preliminaries. In Section 3, we prove the key Theorem \ref{Thm1}. In Section 4, we give the proof of Theorem \ref{Thm2} using Theorem \ref{Thm1}, 



\section{Preliminaries}
In this section, we give some lemmas which will be used to prove Theorems \ref{Thm2} and \ref{Thm1}.

\begin{Lemma}(Dirac \cite{Dirac})\label{Dirac}
Let $G$ be a connected graph with $\delta(G)=\delta$, then
$$p(G)\ge min\left\{|G|, 2\delta+1\right\}.$$
\end{Lemma}

\begin{Lemma}(Chv\'{a}tal and Erd\H{o}s \cite{chvatal})\label{chvatal}
Let $G$ be a graph so that $\alpha(G)\le \kappa(G)+1$, then $G$ has a hamiltonian path.
\end{Lemma}

\begin{Lemma}(Chen et al. \cite{chen})\label{chen0}
Let $G$ be a connected graph of order $n\geq \sum_{i=1}^kn_i$. If $\sigma_3(G)\geq 3\sum_{i=1}^k\lfloor \frac{n_i}{2}\rfloor$, then $G$ contains linear forest $P_{n_1}\cup P_{n_2}\cup \cdots \cup P_{n_k}$.
\end{Lemma}


\begin{Lemma}(Chen et al. \cite{chen})\label{chen}
Let $G$ be a connected graph of order $n$, then the following two results hold.
\begin{enumerate}[\rm(i)]
  \item  If $G$ contains exactly two end blocks, then $p(G)\ge min\{n, \sigma_3(G)\}$.
  \item If $G$ contains only one cut vertex and there is a block that has no hamiltonian path starting at this cut vertex, then $p(G)\ge min\{n, \sigma_3(G)\}$.
\end{enumerate}
\end{Lemma}
\begin{Lemma}(Saito \cite{saito})\label{saito}
Suppose $G$ is a $2$-connected graph of order $n$. Then either $G$ contains a strong dominating cycle or $p(G)\ge min\{n,\sigma_3(G)-1\}$.
\end{Lemma}


\begin{Lemma}\label{claw}
Let $F$ be a linear forest so that $\ell_i\not=3$ for $1\le i\le k$ and $k\ge2$. If $P_{2\delta_F+1}\cup P_{\delta_F}$ is $F$-free, then $F\in \{5P_5, 2P_\ell, P_{\ell+1}\cup P_\ell\}$,  or $F= P_{\ell+2}\cup P_\ell, 2P_{\ell}\cup P_2$ and $\ell$ is odd,
\end{Lemma}
\pf It is obvious that $\ell_1\leq 2\delta_F+1$, and hence $P_{\ell_1}\subseteq P_{2\delta_F+1}$. Suppose $t$ is the maximum number such that $P_{\ell_1}\cup \cdots \cup P_{\ell_t}\subseteq P_{2\delta_F+1}$. Because $F\not\subseteq P_{2\delta_F+1}\cup P_{\delta_F}$, $P_{\ell_{t+1}}\cup\cdots\cup P_{\ell_k}\not\subseteq P_{\delta_F}$. Hence we have
\begin{align}
&\ell_1+\cdots+\ell_{t+1}\ge 2(\delta_F+1)= \sum\limits_{i=1}^k\ell_i-\sum\limits_{i=1}^k\chi_{\mathcal{O}}(\ell_i),\label{eq2}\\
&\ell_{t+1}+\cdots+\ell_k\ge \delta_F+1= \frac{1}{2}\left(\sum\limits_{i=1}^k\ell_i-\sum\limits_{i=1}^k\chi_{\mathcal{O}}(\ell_i)\right),\label{eq3}
\end{align}
By (\ref{eq2}) and (\ref{eq3}), we have
\begin{align}
\ell_{t+1}+3\chi_{\mathcal{O}}(\ell_{t+1}) &\ge \big(\ell_1-3\chi_{\mathcal{O}}(\ell_1)\big)+\cdots+\big(\ell_t-3\chi_{\mathcal{O}}(\ell_t)\big)+\notag\\
&~~~~\big(\ell_{t+2}-3\chi_{\mathcal{O}}(\ell_{t+2})\big)+\cdots+\big(\ell_k-3\chi_{\mathcal{O}}(\ell_k)\big).\label{eq4}
\end{align}
\vskip 2mm

 Assume that $\ell_{t+1}$ is even. If $t+1\geq 3$, then by (\ref{eq4}) we have
 $$\ell_{t+1}\ge(\ell_1-3\chi_{\mathcal{O}}(\ell_1))+(\ell_2-3\chi_{\mathcal{O}}(\ell_2))>\ell_{t+1},$$
or $F=2P_5\cup P_4$ which contradicts (\ref{eq3}). Hence, $t+1=2$.  Note that  (\ref{eq2}) implies $k=t+1=2$  and (\ref{eq3}) implies $\ell_1\le \ell_2+1$.  Let $\ell_2=\ell$ be even, then  $F= 2P_{\ell}, P_{\ell+1}\cup P_{\ell}$.

\vskip 3mm

If $\ell_{t+1}$ is odd, then $\ell_{t+1}\ge 5$. If  $t+1\ge 4$, then we have $t+1\not=k$ by (\ref{eq3}) and
$$\ell_{t+1}+3\ge(\ell_1-3\chi_{\mathcal{O}}\ell_1))+\cdots+(\ell_3-3\chi_{\mathcal{O}}(\ell_3))+\cdots+(\ell_k-3\chi_{\mathcal{O}}(\ell_k))$$
 by (\ref{eq4}), which holds only if $t+1=4$ and $F\in\{5P_5, 4P_5\cup P_2\}$. Because $4P_5\cup P_2\subset P_{2\delta_F+1}\cup P_{\delta_F}$,  
 we have $F=5P_5$. If $t+1\le 3$, then by (\ref{eq2}), we have $k=t+1$ or $F=P_{\ell_1}\cup \cdots\cup P_{\ell_{t+1}}\cup P_2$.  In the former case, by (\ref{eq3}) we can deduce that
 $$\ell_{t+1}+1\ge \left(\ell_1-\chi_{\mathcal{O}}(\ell_1)\right)+\cdots+\left(\ell_t-\chi_{\mathcal{O}}(\ell_t)\right),$$
which implies $t+1=2$ and $\ell_1\le \ell_2+2$. that is,  $F\in \{2P_{\ell_2},P_{\ell_2+1}\cup P_{\ell_2}, P_{\ell_2+2}\cup P_\ell\}$. In the latter case, by (\ref {eq3}) we can deduce that $F\in \{3P_5\cup P_2, 2P_{\ell_2}\cup P_2, P_{\ell_2+1}\cup P_{\ell_2}\cup P_2, P_{\ell_2+2}\cup P_{\ell_2}\cup P_2,P_{\ell_2+3}\cup P_{\ell_2}\cup P_2, P_{\ell_2+4}\cup P_{\ell_2}\cup P_2\}$. It is easy to check that only if $F=2P_{\ell_2}\cup P_2$, $P_{2\delta_F+1}\cup P_{\delta_F}$ is $F$-free. Let $\ell_2=\ell$ be odd. Then we can see that $F\in \{5P_5, 2P_\ell, 2P_{\ell}\cup P_2, P_{\ell+1}\cup P_\ell, P_{\ell+2}\cup P_\ell\}$.


Thus, the result follows by the arguments above.
$\hfill\blacksquare$
\vskip 3mm

\begin{Lemma}\label{claw2}
Let $F$ be a linear forest so that $\ell_i\not=3$ for $1\le i\le k$ and $k\ge2$. If $P_{2\delta_F+1}\cup P_{\delta_F+1}$  is $F$-free, then $F=2P_\ell$ and $\ell$ is odd. Furthermore, $P_{2\delta_F+1}\cup P_{\delta_F+2}$ contains $F$.
\end{Lemma}
\pf Because $P_{2\delta_F+1}\cup P_{\delta_F}$ is a subgraph of $P_{2\delta_F+1}\cup P_{\delta_F+1}$, by Lemma \ref{claw}, it is sufficient to consider $F=5P_5,2P_\ell, P_{\ell+1}\cup P_{\ell}$, $P_{\ell+2}\cup P_{\ell}, 2P_{\ell}\cup P_2$. After an easy check, we find that $P_{2\delta_F+1}\cup P_{\delta_F+1}$ is $F$-free only for $F=2P_{\ell}$ and $\ell$ is odd.
The latter part follows from the fact that $P_{2\delta_F+1}\cup P_{\delta_F+2}$ contains $2P_\ell$. $\hfill\blacksquare$


\begin{Lemma}\label{lem8-1} Let $F$ be a linear forest so that $F\not=kP_5$ and $\ell_i\not=3$ for $1\leq i\leq k$. Suppose $K_{s,t}=(X,Y)\subseteq G$ with $|X|=\delta_F$ and $|Y|\geq \sum_{i=1}^k\ell_i-\delta_F$. Then $G$ has an $F$ if some $\ell_i$ is even and $v_1v_2\in G-(X\cup Y)$ such that $x_1\in N_X(v_1)$, or each $\ell_i$ is odd and $G-(X\cup Y)$ has a path $v_1v_2v_3$ such that $x_1\in N_X(v_1)$ or two independent edges $v_1v_2,v_3v_4$ such that $v_1x_1,v_3x_2\in [\{v_1,v_3\},X] $ and $x_1\not=x_2$.
\end{Lemma}

\pf By the assumption, we let $\ell_i$ be even or $\ell_i\geq 7$ be odd. Partition $(X,Y)$ into $(X_1,Y_1),(X_2,Y_2),...,(X_k,Y_k)$ such that $|X_j|=\lfloor\frac{\ell_j}{2}\rfloor$ and $|Y_j|=\lceil\frac{\ell_j}{2}\rceil$ if $j\not=i$, $|X_i|=\lfloor\frac{\ell_i}{2}\rfloor-1$ and $Y_i=Y\backslash (\cup_{j\not=i} Y_j)$. Clearly, $|Y_i|\geq \lceil\frac{\ell_i}{2}\rceil+1$ and $(X_j,Y_j)$ has a $P_{\ell_j}$ if $j\not=i$. If $\ell_i=2$, then
$v_1v_2$ is a $P_{\ell_i}$. If $\ell_i \geq 4$ be even, we let $x_1\in X_i$. Then $(X_i,Y_i)$ has a $P_{\ell_i-2}$ starting at $x_1$, which together with $v_1v_2$ gives a $P_{\ell_i}$ in $G$. If $\ell_i\geq 7$ is odd, we let $x_1,x_2\in X_i$. Then $(X_i,Y_i)$ has a $P_{\ell_i-3}$ starting at $x_1$ and a $P_{\ell_i-4}$ starting at $x_1$ and ending at $x_2$. Thus, the path $P_{\ell_i-3}$ together with $v_1v_2v_3$ or the path $P_{\ell_i-4}$ together with $v_1v_2,v_3v_4$ gives a $P_{\ell_i}$ in $G$. Therefore, the conclusion holds.$\hfill\blacksquare$

\begin{Lemma}\label{lem8-2} Let $F$ be a linear forest so that $F\not=kP_5$  and $\ell_i\not=3$ for $1\le i\le k$  and $G$ an $F$-free connected graph with $|G|\ge 2|F|$. If $G$ contains a strong dominating path $P$ such that $d(v)\geq \delta_F$ for any $v\in V(G)\backslash V(P)$, then $G\subseteq G_F(n)$.
\end{Lemma}

\pf  Let $P=v_1v_2\cdots v_p$ be a longest path with the property  and $Y=V(G)\backslash V(P)$. By the assumption, $Y\not=\emptyset$.
By the maximality of $P$, $N_Y(v_1)=N_Y(v_p)=\emptyset$  and any vertex of $Y$ has no consecutive neighbors in $P$, which implies $|P|\ge 2\delta_F+1$. If $|P|= 2\delta_F+1$, then  $N(y)=\{v_2,v_4,\ldots,v_{p-1}\}$ for any $y\in Y$. Let $X=N(y)$. Obviously, $|X|=\delta_F$, $|X|+|Y|\ge |F|$ and $G[X,Y]$ is a complete bipartite graph. By Lemma \ref{lem8-1}, $P-X$ contains at most one edge, according to $F$ consists of only odd paths or not. This implies $G\subseteq G_F(n)$. Hence we may assume  $|P|\ge 2\delta_F+2$.

We proceed on  by induction on $k$.

If $k=1$, then since $\ell_1-1\geq |P|\ge 2\delta_F+2\geq 2\lfloor \frac{\ell_1}{2}\rfloor$, we deduce that $|P|=\ell_1-1$ and $\ell_1\geq 7$ is odd. In this case, $d(y)=\delta_F=d$ for any $y\in Y$. 
Let $y\in Y$ be given and $N(y)=X=\{x_1,x_2,...,x_d\}$, where the subscripts occur along $P$ from $v_1$ to $v_p$. Set $Q_i=x_i\ora{P}x_{i+1}$ for $0\leq i\leq d$, where $x_0=v_1$ and $x_{d+1}=v_p$. Suppose there is some $y'\in Y\backslash \{y\}$ such that $N(y')\not=X$. Let $a\in N_{Q_i}(y')\backslash X$.
Obviously, $|Q_i|=3$ or 4. If $i=0$ or $i=d$, then $|Q_i|=3$ and $|Q_j|=3$ for $1\leq j\leq d-1$. By symmetry of $Q_0$ and $Q_d$, we assume that $i=0$. By the maximality of $P$,  $x_1\notin N(y')$ and $N(y')\cap \{x_1^+,x_d\}\not=\emptyset$. Thus, either $x_0ay'x_1^+x_1yx_2\ora{P}x_{d+1}$ is a $P_{\ell_1+1}$, or $x_0ay'x_d\ola{P}x_1y$ is a $P_{\ell_1}$, a contradiction. Hence we have $1\leq i\leq d-1$. If $|Q_i|=3$, then by the maximality of $P$, $x_i,x_{i+1}\notin N(y')$ and $N(y')\cap \{x_i^-,x_{i+1}^+\}\not=\emptyset$. Thus, $x_0\ora{P}x_i^-y'ax_iyx_{i+1}\ora{P}x_{d+1}$  or $x_0\ora{P}x_iyx_{i+1}ay'x_{i+1}^+\ora{P}x_{d+1}$ is a  $P_{\ell_1+1}$, a contradiction. If $|Q_i|=4$, then $a=x_i^+$ or $x_{i+1}^-$. By symmetry, we assume that $a=x_i^+$.
By the maximality of $P$, $x_i,x_{i+1}^-\notin N(y')$. If $x_i^-\in N(y')$, then $x_0\ora{P}x_i^-y'ax_iyx_{i+1}\ora{P}x_{d+1}$ is a  $P_{\ell_1}$ and so $x_i^-\notin N(y')$, which implies $x_{i+1}\in N(y')$. Thus, $x_0\ora{P}x_iyx_d\ola{P}x_{i+1}y'ax_{i+1}^-$ is a $P_{\ell_1}$, a contradiction. Therefore, $N(y')=X$ for any $y'\in Y$. By Lemma \ref{lem8-1}, $P-X$ has at most one edge, and hence we have $G\subseteq G_F(n)$.

Assume that $k\geq 2$ and the conclusion holds  for $k-1$.  Consider $F=P_{\ell_1}\cup\cdots\cup P_{\ell_k}$. 
Let $F'=F-P_q$, where $q=5$ if some $\ell_i=5$, and $q=\ell_k$ otherwise. Set $P'=v_1v_2\cdots v_q$, $P''=P-V(P')$  and $G'=G- V(P')$.
For any $y\in Y$, noting that $d_{P'}(y)\le \lfloor\frac{q}{2}\rfloor$, we have  $d_{P''}(y)\ge \delta_{F'}$. Since
 $|P|\geq 2\delta_F+2$, we have $|P''|\geq 2\delta_F+2-|P'|\ge 2\delta_{F'}+1$. Clearly,   $G'$ is $F'$-free and $|G'|=n-q\ge 2|F'|$,  and $P''$ is a strong dominating path of $G'$ such that $d_{P''}(y)\ge \delta_{F'}$ for any $y\in Y$. By induction hypothesis, $G'\subseteq G_{F'}(n-q)$. Let $V(G')=U\cup V$ with $|U|=\delta_{F'}$ and $V=V(G')\backslash U$ such that  $U=V(K_{\delta_{F'}})$ in $G_{F'}(n-q)$.
 By the definition of $G_{F'}(n-q)$, we can see that any path of order at least $2\delta_{F'}+1$ in $G'$ must contain all vertices of $U$.
 Thus, we have $U\subseteq V(P'')$ and so  $Y\subseteq V$. Since $G'[V]$ has at most one edge, $Y$ is an independent set and $d(y)\geq \delta_{F'}$ for any $y\in Y$, we can deduce that  all vertices in $Y$ have the same neighborhoods on $P''$ except at most one vertex $y'\in Y$, and $d_{P''}(y)=\delta_{F'}$ for any $y\in Y\backslash \{y'\}$.
Let $Q''=v_1v_2\cdots v_{p-q}$ and $Q'=P-Q''$. By the symmetry of $P''$ and $Q''$, all vertices in $Y$, except at most one vertex $y''\in Y$, have the same neighborhoods on $Q''$, and $d_{Q''}(y)=\delta_{F'}$ for any $y\in Y\backslash \{y''\}$.



If $|P'|=\ell_k$ is even, then since $d_{P''}(y)=\delta_{F'}$ and $d(y)\ge \delta_F$ for any $y\in Y\backslash \{y'\}$, we can  get $d_{P'}(y)=\frac{\ell_k}{2}$, and hence $N_{P'}(y)=\{v_2,\ldots, v_{\ell_k}\}$.  Thus,  $N_P(y_1)=N_P(y_2)$ for any $y_1,y_2\in Y\backslash\{y'\}$.
If $|P'|=\ell_k$ is odd or $|P'|=5$, then
noting that $|P''|\ge2\delta_{F'}+1$, we have
 $$|Q''|=|P''|\ge 2\delta_{F'}+1=\sum_{i=1}^{k}\left(\ell_i-\chi_B(\ell_i)\right)-(q-\chi_B(q))-1\ge |P’|,$$
unless $F=2P_{\ell_1}$ and $\ell_1$ is odd, or $F\in\{P_5\cup P_2, P_5\cup P_4, P_5\cup 2P_2\}$. If $|Q''|\geq |P'|$, then  $P'\subseteq Q''$. Thus,  $N_{P'}(y_1)=N_{P'}(y_2)$ and hence $N_P(y_1)=N_P(y_2)$ for any $y_1,y_2 \in Y\backslash\{y',y''\}$.
If $F=2P_{\ell_1}$ and $\ell_1$ is odd, then
$2\delta_{F'}+2\ge|Q''|=|P''|\ge2\delta_{F'}+1$. Let $y_1,y_2 \in Y\backslash\{y',y''\}$, then $d_{P''}(y_i)=d_{Q''}(y_i)=\lfloor\ell_1/2\rfloor-1$ for $i=1,2$ by the arguments above. Since $d(y_i)\geq \delta_F=2\lfloor \frac{\ell_1}{2}\rfloor-1$ and $|P-Q''-P''|\leq 2$, $y_i$ has exactly one neighbor in $P-Q''-P''$.
If $|P''|=2\delta_{F'}+1=\ell_1-2$, then $|P|=2\ell_1-2$ and $V(P-Q''-P'')=\{v_{\ell_1-1},v_{\ell_1}\}$.  If $y_1,y_2$ have different neighbors in $\{v_{\ell_1-1},v_{\ell_1}\}$, it is clear that $G$ contains $2P_{\ell_1}$,
and so $y_1,y_2$ have the same neighbor in $\{v_{\ell_1-1},v_{\ell_1}\}$. Thus we have $N_P(y_1)=N_P(y_2)$.
If $|P''|=2\delta_{F'}+2=\ell_1-1$, then $|P|=2\ell_1-1$ and $V(P-Q''-P'')=\{v_{\ell_1}\}$, and so $v_{\ell_1}\in N_P(y_1)\cap N_P(y_2)$, which implies that $N_P(y_1)=N_P(y_2)$. If $F\in\{P_5\cup P_2, P_5\cup P_4, P_5\cup 2P_2\}$, then since $|P|\geq 2\delta_F+2$ and $d_P(y)\geq \delta_F$ for any $y\in Y$, it is easy to check  that $G$ contains $F$, a contradiction.

By the arguments above, we have $N_P(y_1)=N_P(y_2)$ for any $y_1,y_2 \in Y\backslash\{y',y''\}$.

Let $Y'=Y\backslash\{y',y''\}$, $y\in Y'$ be given and $X=N(y)$, then $|X|=\delta_F$. Obviously, $|Y'|\ge |F|-|X|$  and the edges between $X$ and $Y'$ form a complete bipartite graph in $G$. Since $G$ is $F$-free, by applying Lemma \ref{lem8-1} on $G-\{y',y''\}$,  we have $G-\{y',y''\}\subseteq G_F(n)$. Add $y'$ and $y''$ (if they exist) to $G-\{y',y''\}$, then since $G$ is $F$-free, it is not difficult to show  $G\subseteq G_F(n)$. $\hfill\blacksquare$

\vskip 3mm
To discuss the structure of an $F$-free graph for $F=kP_5$, we need to define a connected graph $G_F(n,i)$ of order $n$ as follows. Add $i$ new independent edges $e_1$,...,$e_i$ to $K_{\delta_F}+E_{n-\delta_F-2i}$ and then connect $V(\cup_{j=1}^ie_j)$ and $V(K_{\delta_F})$  such that $[V(e_j),V(K_{\delta_F})]$ contains two independent edges for each $e_j$  except at most two of $e_1,e_2,...,e_i$, and any vertex in $V(K_{\delta_F})$ has neighbors in at most one $e_j$, illustrated in Figure 2.

\begin{picture}(70,40)(-90,70)
\centering
{\includegraphics[width=2.8in]{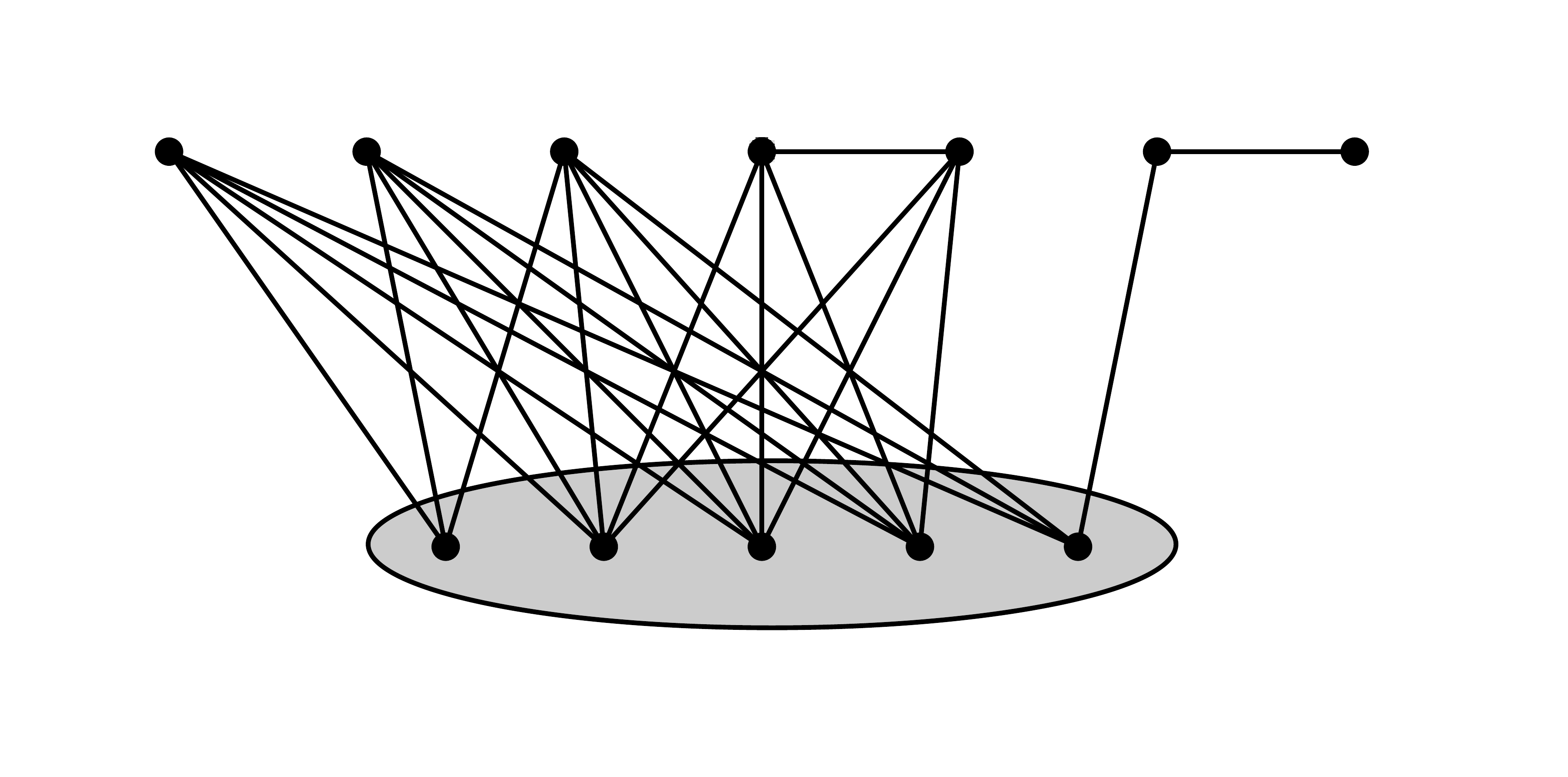}}
\put(-110,-5){\makebox(3,3){Figure 2. $G_F(n,i)$ (oval denotes a clique $K_{\delta_F}$)}}
\put(-42,90){\makebox(4,4){$e_1$}}
\put(-67,90){\makebox(4,4){$\cdots$}}
\put(-94,90){\makebox(4,4){$e_i$}}
\end{picture}
\\
\vskip 20mm

\begin{Lemma}\label{lem8-3}
 Let $F=kP_5$ with $k\geq 2$ and $G$ be an $F$-free connected graph with $|G|\ge 2|F|$. If $G$ contains a strong dominating path $P$ such that $d(v)\geq \delta_F$ for any $v\in V(G)\backslash V(P)$, then $G\subseteq G_F(n,i)$, where $i\le k$.
\end{Lemma}
\pf  Let $P=v_1v_2\cdots v_p$ be a longest path with this property and $Y=V(G-P)$. By the maximality of $P$, $N_Y(v_1)=N_Y(v_p)=\emptyset$  and any vertex of $Y$ has no consecutive neighbors in $P$, and so $|P|\ge 2\delta_F+1$. If $|P|= 2\delta_F+1$, then  $N(y)=\{v_2,v_4,\ldots,v_{p-1}\}$ for any $y\in Y$. This implies that $G\subseteq K_{\delta_F}+E_{n-\delta_F}$ and hence the conclusion holds. Thus, we may assume  $|P|\ge 2\delta_F+2$.

We proceed on by induction on $k$.

Assume that  $k=2$. Then $\delta_F=3$ and $8\le p\le 9$.  Let $y,y'\in Y$. We first show that $N(y)=N(y')$. Suppose $p=8$. If $\{v_2,v_7\}\subset N(y)$ for some $y$,  then by the symmetry we may assume $N(y)=\{v_2,v_4,v_7\}$.  If $N(y’)\cap \{v_3,v_5,v_6\}\not=\emptyset$, then we can find $2P_5$, a contradiction. Hence $N(y’)=N(y)$.  If $\{v_2,v_7\}\not\subset N(y)$ for any $y$, then either $N(y)=\{v_2,v_4,v_6\}$ or $N(y)=\{v_3,v_5,v_7\}$. If $N(y)=\{v_2,v_4,v_6\}$ and $N(y')=\{v_3,v_5,v_7\}$, then we can find $2P_5$ easily, also a contradiction. Hence $N(y)=N(y')$. Assume that $p=9$. Since $G$ is $F$-free, $N(y)\cap \{v_4,v_6\}=\emptyset$ and $\{v_2,v_8\}\not\subseteq N(y)$ for any $y\in Y$. Thus, by symmetry we may assume  $N(y)=\{v_3,v_5,v_7\}$ or $N(y)=\{v_3,v_5,v_8\}$. If $N(y')\not=N(y)$, then we can find $2P_5$ in $G[V(P)\cup \{y,y'\}]$, and hence $N(y)=N(y')$.  Therefore, we always have $N(y)=N(y')$.  Now let $X=N(y)$, then $|X|=\delta_F$, $|X|+|Y|\ge 10$ and $[X,Y]$ is a complete bipartite. Consider the vertices $P-X$ and note that $G$ has no $2P_5$, we have $G\subset G_F(n,i)$ with $i\leq 2$.

Assume that $k\geq 3$ and the conclusion holds  for $k-1$. Let  $F’=(k-1)P_5$, $P’=v_1v_2\cdots v_5$ and $P''=P-P’$. Since $d_{P'}(y)\leq 2$ for any $y\in Y$, $P''$ is a strong dominating path of $G-P’$ and $d_{P''}(v)\ge \delta_{F’}$. Clearly, $G-P’$ is $F’$-free. By induction hypothesis, $G-P’\subset G_{F’}(n-5,i)$ with $i\leq k-1$. By the definition of $G_{F’}(n,i)$, any strong dominating path of $G-P'$ must contain all vertices of the $i$ independent edges except at most one vertex $y'$ and all vertices in the clique $K_{\delta_{F’}}$, and so all vertices in $Y-\{y’\}$ have the same neighbors on $P''$. Let $Q''=v_1v_2\cdots v_{p-5}$, then $P’\subset Q''$. By symmetry of $Q''$ and $P''$, all vertices in $Y-\{y''\}$ have the same neighbors on $Q''$ and hence on $P’$. Let $Y'=Y\backslash\{y',y''\}$. By the arguments above, we have $N_P(y_1)=N_P(y_2)$ for any $y_1,y_2 \in Y'$. Let $X=N(y)$ for some $y\in Y'$, then $|X|=\delta_F$, $|X|+|Y'|\ge 5k$ and $[X,Y']$ is a complete bipartite. Consider the vertices in $P-X$ and note that $G$ is $F$-free, we have $G\subset G_F(n,i)$, $i\le k$. $\hfill\blacksquare$

\section{$F$-free graphs with large minimum degree}

In this section, we will give the proof of Theorem \ref{Thm1}, which characterizes the structure of an $F$-free graph $G$ with $\delta(G)\ge \delta_F$.
\vskip 2mm
\noindent{\bf Proof of Theorem \ref{Thm1}.} Suppose that $G$ is 2-connected. If $G$ has a strong dominating path, then $G\subseteq G_F(n)$ for $F\not=kP_5$ by Lemma \ref{lem8-2}, and $G\subseteq G_F(n,i)$ for $F=kP_5$ by Lemma \ref{lem8-3}. In the latter case, we have $i\leq 1$ because $\delta(G)\geq \delta_F$ and hence $G\subseteq G_F(n,1)\subseteq G_F(n)$.
Therefore, we may assume that $G$ has no dominating path. 

By Lemma \ref{saito}, $p(G)\ge min\{n,\sigma_3(G)-1\}$. Since
$G$ is $F$-free, we have
\[\begin{split}
|F|>\sigma_3(G)-1&\ge\frac{3}{2}\sum\limits_{i=1}^{k}(\ell_i-\chi_{\mathcal{O}}(\ell_i))-4=\sum_{i=1}^k\ell_i+\frac{1}{2}\sum\limits_{i=1}^{k}(\ell_i-3\chi_{\mathcal{O}}(\ell_i))-4,\\
\end{split}\]
which implies 
$$F\in \{3P_5,2P_5\cup P_2, P_5\cup2P_2,3P_2, P_7\cup P_5, P_7\cup P_2,P_5\cup P_4,P_4\cup P_2,2P_5, P_5\cup P_2,2P_2\}.$$

Let $P=v_1v_2\cdots v_p$ be a longest path in $G$. By Lemma \ref{saito},  $p\geq 3\delta_F-1\geq |F|-2$. Assume $G-P$ has a component $H$ with $|H|\geq 2$. Since $G$ is 2-connected, there exist two independent edges $x_1v_i,x_tv_j$ with $x_1,x_t\in V(H)$ and $i<j$. Let $P'=x_1x_2\cdots x_t$ be a  path connecting $x_1$ and $x_t$ in $H$. Obviously, $t\geq 2$. By the maximality of $P$, we have 
\begin{equation}\label{eq5}
t+1\leq i\leq j-t-1<j\leq p-t.
\end{equation}
If $F\in \{2P_5\cup P_2,3P_2, P_5\cup P_2,2P_2, P_5\cup 2P_2,P_4\cup P_2, P_7\cup P_2\}$, then since $p\geq |F|-2$, we can see that $G$ contains $F$, a contradiction. 
Hence, $F\in \{P_5\cup P_4,P_5\cup P_7,2P_5,3P_5\}$. 
If $i=3$ or  $j=p-2$, then it is easy to see that $P\cup P'\cup \{x_1v_i,x_tv_j\}$ contains $F$, and hence $i\ge 4$ and $j\le p-3$. This implies $p\ge 10$ and so $F\in \{P_5\cup P_7,3P_5\}$. If $F=P_5\cup P_7$, then $\delta_F=4$ and $p=11$. By  (\ref{eq5}), $i=4$ or $j=8$, and hence $P\cup P'\cup \{x_1v_i,x_tv_j\}$ contains $P_5\cup P_7$, a contradiction. Thus, $F=3P_5$. In this case, $\delta_F=5$ and $p=14$. For any edge $x'x''\in G-P$, since $G$ has no $3P_5$, we have $N_P(x')\cup N_P(x'')\subseteq \{v_5,v_{10}\}$. 
Noting that $\delta(G)\geq 5$, we get $\delta(H)\geq 3$. Because $H$ contains no $P_5$,  $H=K_4$ and $N_P(x)=\{v_5,v_{10}\}$ for any $x\in V(H)$. By symmetry, we have $G[\{v_6,v_7,v_8,v_9\}]=K_4$. Thus, $G-V(P)-V(H)$ has no isolated vertices because such a vertex has at most 4 neighbors in $P$ which contradicts $\delta(G)\geq 5$.  Note that $\delta(G)\geq 5$ and $G$ has no $3P_5$, we  can see $G[\{v_1,v_2,v_3,v_4\}]=G[\{v_{11},v_{12},v_{13},v_{14}\}]=K_4$, and $G=K_2+tK_4$ or $E_2+tK_4$.

%

\vskip 4mm

Next we assume that $G$ is not 2-connected. We distinguish the following two cases separately.
\vskip 2mm
\noindent{\bf Case 1.} $G$ contains exactly two end blocks or  $G$ contains only one cut vertex.
\vskip 2mm

If $G$ contains exactly two end blocks, or $G$ has only one cut vertex and some  block contains no hamiltonian path starting at the cut vertex, then by Lemma \ref{chen}, $p(G)\ge  min\{n,\sigma_3(G)\}$. Since $G$ is $F$-free, we have
$$|F|>\sigma_3(G)=\sum_{i=1}^k\ell_i+\frac{1}{2}\sum\limits_{i=1}^{k}(\ell_i-3\chi_{\mathcal{O}}(\ell_i))-3,$$
which implies $F\in\{2P_5,P_5\cup P_2,2P_2\}$. Let $P=v_1v_2\cdots v_p$ be a longest path in $G$, we have $p=\sigma_3(G)=|F|-1$.  If $F=P_5\cup P_2$, then $p=6$ and $G-P$ has no edges, which implies $d_P(v)\geq \delta_F=2$ for any $v\in G-P$.  In this case, it is easy to find an $F$ in $G$, a contradiction. If $P=2P_2$, then $G$ is a star and the result holds trivially. If $F=2P_5$, then $\delta_F=3$ and $p=9$. If $P$ is a strong dominating path, then $G\subset G_F(n,i)$ with $i\le 2$ by Lemma \ref{lem8-3}. Noting that $\delta(G)\geq 3$, we find $G$ is 2-connected, a contradiction.  Hence $P$ is not a strong dominating path. Since $G$ has no $2P_5$, $N_P(x)\cup N_P(x')\subseteq \{x_5\}$ for any $xx'\in G-P$.
Let $H$ be any component of $G-P$ with $|H|\ge2$. Then $N_P(x)\subset\{v_5\}$ for any $x\in V(H)$ and so $\delta(H)\geq 2$. Since $H$ is $P_5$-free, we have $|H|=3,4$ by Lemma \ref{Dirac}.
In addition, if $G-P$ has an isolated vertices $x$, then $\{v_2,v_8\}\not\subset N(v)$ for otherwise $G$ has $2P_5$,  and so we can assume $N(v)=\{v_3,v_5,v_7\}$ or $\{v_3,v_5,v_8\}$ by symmetry. In both case, $G$ contains $2P_5$, a contradiction. Thus, $G-P$ has no isolated  vertices. If there is an edge between $\{v_1,v_2,v_3,v_4\}$ and $\{v_6,v_7,v_8,v_9\}$, then since $G$ has no $2P_5$, the edge must be $v_4v_6$. In this case, $v_4,v_5,v_6$ are cut vertices which contradicts our assumption. So there is no edge between $\{v_1,v_2,v_3,v_4\}$ and $\{v_6,v_7,v_8,v_9\}$. Therefore, $v_5$ is the only cut vertex, $G-v_5$ has at least three components and each component has a hamiltonian path starting at $v_5$, a contradiction.

 \vskip 2mm
Next we may assume that $v$ is the only cut vertex of $G$, $B_1,\ldots, B_t$ are all blocks and $t\ge 3$, and each  $B_i$ contains a hamiltonian path starting at $v$. Obviously, $|B_i|\ge \delta_F+1$ for $1\le i\le t$.  Assume that $|B_1|\geq \cdots \geq |B_t|$ and let $Q_i$ be a hamiltonian path starting at $v$ in $B_i$ for $1\leq i\leq t$. Because $|Q_i|\ge \delta_F+1$, $|Q_1-v|\geq \delta_F$ and $Q_2\cup Q_3$ is a path of order at least $2\delta_F+1$. Since $G$ is $F$-free,  by Lemma \ref{claw2}, we have $|Q_1-v|\leq \delta_F+1$.

If $|Q_1-v|=\delta_F$, then by Lemma \ref{claw}, $F\in\{5P_5, 2P_{\ell},P_{\ell+1}\cup P_{\ell}\}$, or $F\in\{P_{\ell+2}\cup P_{\ell}, 2P_\ell\cup P_2\}$ with $\ell$ odd. In this case, noting that $\delta(G)\geq \delta_F$, we can see that $G=K_1+tK_{\delta_F}$ with $t\geq 4$(because $n\ge 2|F|$). Note that $G$ contains $F$ if $F=5P_5$, hence $F\in\{2P_{\ell},P_{\ell+1}\cup P_{\ell}\}$, or $F\in\{P_{\ell+2}\cup P_{\ell}, 2P_\ell\cup P_2\}$ with $\ell$ odd. If $|Q_1-v|=\delta_F+1$, then by Lemma \ref{claw2}  and the fact that $P_{2\delta_F+1}\cup P_{\delta_F+1}$  has no
$F$, we have $F=2P_\ell$ and $\ell$ is odd. In this case, it is not difficult to deduce that $G\subseteq K_1+(t_1K_{\ell-1}\cup t_2K_{\ell-2})$.

\vskip 2mm
\noindent{\bf Case 2.} $G$ contains at least three  end blocks and two cut vertices.
\vskip 2mm

By the assumption, we may assume that $B_1,B_2,B_3$ are three end blocks and $b_i\in B_i$ is the only cut vertex of $B_i$ with $b_1\not=b_2$. Let $Q_i$ be the longest path in $B_i$ starting at $b_i$ for $i=1,2,3$,  $P$ a longest path connecting $b_1,b_2$ and $P'$ a path with endpoints $b_3$ and $v$ such that $V(P)\cap V(P')=\{v\}$. Because $\delta(G)\ge \delta_F$, we have $|Q_i|\ge \delta_F+1$. Thus, it is not difficult to see that the subgraph consisting of the paths $Q_1,Q_2,Q_3,P,P'$ contains
$P_{2\delta_F+1}\cup P_{\delta_F+1}$ which is $F$-free. By Lemma \ref{claw2}, $F=2P_{\ell}$ and $\ell$ is odd.

Let $S=V(B_1)\cup V(B_2)\cup V(P)$. It is clear that $N_S(u)\subseteq V(P)$ for any $u\in V(G)\backslash S$.

Obviously, $\delta(G)\geq \delta_F=\ell-2$. Since $G$ is $F$-free, we have $|Q_1\cup P\cup Q_2|\le 2\ell-1$, and hence $|P|\leq 3$. Note that $|P|\geq 2$, we have $|P|=2$ or 3.

If $|P|=2$, then $P=b_1b_2$. For any $u\in V(G)\backslash S$, by the maximality of $P$, $d_P(u)\leq 1$, which means $N_S(u)\subseteq \{b_1\}$ or $N_S(u)\subseteq \{b_2\}$. If $u_1,u_2\in V(G)\backslash S$ with $u_1b_1,u_2b_2\in E(G)$, then $G$ contains $2P_\ell$, and hence we may assume that  $N_S(u)\subseteq \{b_1\}$. In this case, $G-b_1$  contains at least $3$ components and each one is $P_\ell$-free for otherwise $G$ has $2P_\ell$. Since $\delta(G)\geq \ell-2$,  each component contains $\ell-2$ or $\ell-1$ vertices, $B_2=K_{\ell-1}$ and $G-V(B_2)\subseteq K_1+(s_1K_{\ell-1}\cup s_2K_{\ell-2})$. Of course, $G\subseteq K_1+(t_1K_{\ell-1}\cup t_2K_{\ell-2})$.

If $|P|=3$, then $N(u)\cap \{b_1,b_2\}=\emptyset$ for otherwise $G$ has $2P_\ell$. Thus, $P=b_1vb_2$ and  $N_S(u)\subseteq \{v\}$ for any $u\in V(G)\backslash S$. If $|B_1|\geq \ell$, then since $\delta(B_1)\geq \ell-3$, by Lemma \ref{Dirac} we have $p(B_1)\geq min\{|B_1|, 2(\ell-3)+1\}\geq \ell$,  that is, $B_1$ contains a $P_\ell$. Note that $Q_2\cup \{v\}$ also contains a $P_\ell$, a contradiction. Thus, $|B_1|=\ell-1$. By symmetry, $|B_2|=\ell-1$. Since $\delta(G)\geq \ell-1$, we have $B_1=B_2=K_{\ell-1}$. By the assumption, $V(G)\backslash S\not=\emptyset$. Let $G'=G-B_1-B_2$, then $G'$ is connected. By the arguments before, $d_{G'}(v)\geq \delta_F-2$ and $d_{G'}(u)\geq \delta_F$ for any $u\in V(G')\backslash \{v\}$.

If $b_1b_2\in E(G)$, then $Q_1\cup \{b_1b_2\}\cup Q_2$ contains a $P_{\ell}$. If $G’-v$ contains only one component, then $G’$ contains a path of order at least
                $$min\{|G-B_1-B_2|, 2\delta_F-1\}\ge min\{|G-B_1-B_2|, \ell\}.$$
Because $G'$ is $P_\ell$-free, we have $\ell-1\ge |G'|\ge\delta_F+1=\ell-1$, and so $G'=K_{\ell-1}$. In this case, $G$ consists of $B_1,B_2,B_3$ and $b_1b_2b_3$ is a triangle, which is impossible since $|G|=3(\ell-1)<2|F|$. If $G’-v$ contains two components, we can find a $P_\ell$ easily, a contradiction.

If $b_1b_2\not\in E(G)$, then $G-v$  contains at least $3$ components including $B_1,B_2$, and each one is $P_\ell$-free, for otherwise $G$ has  $2P_\ell$. Since $\delta(G)\geq \ell-2$,  each component of $G-v$ contains $\ell-1$ or $\ell-2$ vertices. This is to say $G'\subseteq K_1+(s_1K_{\ell-1}\cup s_2K_{\ell-2})$, and so $G\subseteq K_1+(t_1K_{\ell-1}\cup t_2K_{\ell-2})$.
\vskip 2mm

\
The proof of Theorem \ref{Thm1} is complete. $\hfill\blacksquare$

\section{Generalized Tur\'{a}n number}

In this section, we determine the generalized Tur\'{a}n number $ex(n,K_s,F)$. We will use a strategy to reduce the problem to graphs with large minimum degree. For a fixed $F$, we define a process of $G$, called $\left(\delta_F-1\right)$-$disintegration$ as follows: for $j=n$, let $G^n=G$, and for $j<n$, let $G^j$ be obtained from $G^{j+1}$ by deleting a vertex of degree less than $\delta_F$ in $G^{j+1}$ if such a vertex exists. The process terminates at $G^m$ when $\delta (G^m)\ge \delta_F$ or $m=0$. We call $G^m$ the $\delta_F$-$core$ of $G$.

Before starting to prove Theorem \ref{Thm2}, as a preparation, we first discuss the properties of the $\delta_F$-core of $F$-free graphs for some special $F$.

\begin{Claim}\label{c1}
Let $F$ be a linear forest and $G$ an $F$-free graph with $\delta(G)\ge \delta_F$ and components $G_1,G_2,\ldots,G_c$, where $|G_1|\geq \cdots\geq |G_c|$. Then  $G=G_1\cup (c-1)K_{\ell-1}$ if $F=2P_\ell$ and $\ell$ is odd, and $|G_i|\le |F|-1$ for all $1\le i\le c$ otherwise.
\end{Claim}
\pf If $F=2P_\ell$ and $\ell$ is odd, then $\delta_F=\ell-2$. Since $p(G_i)\ge min\{|G_i|, 2\ell-3\}$ by Lemma \ref{Dirac}, we have $|G_i|\leq \ell-1$  which implies $G_i=K_{\ell-1}$ for $2\leq i\leq c$, and hence the result follows.


If $F\not= 2P_\ell$ with $\ell$ odd, we suppose $|G_1|\ge |F|$. Since $G_1$ contains no hamiltonian path, then $\alpha(G)\ge 3$ by Lemma \ref{chvatal}.
Obviously, 
$$\sigma_3(G_1)\ge 3\delta_F\ge 3(\delta_{F'}+1),$$
where $F'=F-P_{\ell_k}$. By Lemma \ref{chen0},  we have $F'\subseteq G_1$. Since $G_2$ contains $P_{\ell_k}$, we can see that $G$ contains $F$,  a contradiction.
$\hfill\blacksquare$

\begin{Claim}\label{c2} Let $F\in \{3P_5,2P_{\ell},P_\ell\cup P_{\ell+1},P_\ell\cup P_{\ell+2}, 2P_{\ell}\cup P_2\}$ and $s\le \delta_F+1$.  \\
If $F\not=2P_\ell$ with $\ell$ odd or $F\not=3P_5$, then
$$\mathcal{N}_s(G_F(t(\ell-1)+1))\ge\mathcal{N}_s(K_1+tK_{\ell-1}).$$
If $F=3P_5$, then
$$\mathcal{N}_s(G_F(4t+2))>\mathcal{N}_s(K_2+tK_4).$$
If $F=2P_\ell$ with $\ell$ odd, $2\le s\le \ell-2$ and  $m\ge 5{\binom{|F|-1}{s}}/{\binom{\delta_F}{s-1}}$, then
$$\mathcal{N}_s(G_F(m))=max\{\mathcal{N}_s(G)~|~G~is~F\text{-}free,~ |G|=m~and~\delta(G)\geq \delta_F\}.$$

%
\end{Claim}
\pf If $F\not=2P_\ell$ with $\ell$ odd or $F\not=3P_5$, then $\delta_F=\ell-1$. Thus we have

\[\begin{split}
\mathcal{N}_s(G_F(t(\ell-1)+1))-\mathcal{N}_s(K_1+tK_{\ell-1})&\ge \binom{\ell-1}{s}+(m-\ell+1)\binom{\ell-1}{s-1}-t\binom{\ell}{s}\\
&=(t-1)\left(\ell-1-\frac{\ell}{s}\right)\binom{\ell-1}{s-1}\ge 0,
\end{split}\]
that is, $\mathcal{N}_s(G_F(t(\ell-1)+1))\ge\mathcal{N}_s(K_1+tK_{\ell-1})$.
\vskip 2mm

If $F=3P_5$, then $\delta_F=5$. Thus we have 
\[\begin{split}
\mathcal{N}_s(G_F(4t+2))-\mathcal{N}_s(K_2+tK_{4})&\ge \binom{5}{s}+(4t-3)\binom{5}{s-1}+\binom{5}{s-2}-t\binom{6}{s}\\
&=\big((2s/3-1)(t-1)-s/2\big)\binom{6}{s}+\binom{5}{s-2}> 0.
\end{split}\]

\vskip 2mm
If $F=2P_\ell$ with $\ell$ odd, then $\delta_F=\ell-2$.
Since $m\ge 5{\binom{|F|-1}{s}}/{\binom{\delta_F}{s-1}}$ and $G$ is $F$-free, by Theorem \ref{Thm1} and Claim \ref{c1}, we have $G\in \{G_F(m), K_1+(t_1K_{\ell-1}\cup t_2K_{\ell-2}), G_1\cup tK_{\ell-1}\}$, where $G_1$ is the largest component of $G$ with $|G_1|=n_1$ if $G$ is disconnected.
We first show $\mathcal{N}_s\big(G_F(m)\big)\ge \mathcal{N}_s\big(K_1+(t_1K_{\ell-1}\cup t_2K_{\ell-2})\big)$ when $m=1+t_1(\ell-1)+t_2(\ell-2)$. Noting that $s\geq 2$ and $m\ge 5{\binom{|F|-1}{s}}/{\binom{\delta_F}{s-1}}$ implying $t_1+t_2\geq 5$, we have
\[\begin{split}
&\mathcal{N}_s\big(G_F(m)\big)-\mathcal{N}_s\big(K_1+(t_1K_{\ell-1}\cup t_2K_{\ell-2})\big)\\
=&\binom{\ell-2}{s-2}+(m-\ell+2)\binom{\ell-2}{s-1}+\binom{\ell-2}{s}-t_2\binom{\ell-1}{s}-t_1\binom{\ell}{s}\\
=& \binom{\ell-2}{s-2}+\biggl(t_1(\ell-1)+(t_2-1)(\ell-2)+1\biggr)\binom{\ell-2}{s-1}+\frac{\ell-s-1}{s}\binom{\ell-2}{s-1}\\
&~~~~~~~~~~~-t_1\frac{\ell(\ell-1)}{s(\ell-s)}\binom{\ell-2}{s-1}-t_2\frac{\ell-1}{s}\binom{\ell-2}{s-1}\\
=&\binom{\ell-2}{s-2}+\biggl(t_1(\ell-1)\left(1-\frac{\ell}{s(\ell-s)}\right)+(t_2-1)\left(\ell-2-\frac{\ell-1}{s}\right)\biggr)\binom{\ell-2}{s-1}\\
\geq&\binom{\ell-2}{s-2}+\biggl(t_1(\ell-1)\frac{\ell-4}{2(\ell-2)}+(t_2-1)\frac{\ell-3}{2}\biggr)\binom{\ell-2}{s-1}>0.
\end{split}\]

Next we show $\mathcal{N}_s(G_F(m))\ge\mathcal{N}_s(G_1\cup tK_{\ell-1})$ when $m=n_1+t(\ell-1)$.  If $n_1\ge 2|F|$, then $G_1\subseteq G_F(n_1)$ by Lemma \ref{lem8-2}. Thus we have
\[\begin{split}
\mathcal{N}_s\left(G_F(m)\right)-\mathcal{N}_s(G_1\cup tK_{\ell-1})
\geq&\mathcal{N}_s(G_F(m))-\mathcal{N}_s(G_F(n_1))-t\binom{\ell-1}{s}\\
=&t(\ell-1)\binom{\ell-2}{s-1}-t\binom{\ell-1}{s}>0.
\end{split}\]
If $\ell-1\le n_1\le 2|F|-1$, then noting that $G_1$ is $P_{2\ell}$-free,  by Theorem \ref{Luo} we have
$$ \mathcal{N}_s(G_1)\le \frac{n_1}{2\ell-1}\binom{2\ell-1}{s}.$$
Thus, note that $m\ge 5{\binom{|F|-1}{s}}/{\binom{\delta_F}{s-1}}$, that is, $m\ge 5{\binom{2\ell-1}{s}}/{\binom{\ell-2}{s-1}}$, we have

\[\begin{split}
&\mathcal{N}_s(G_F(m))-\mathcal{N}_s(G_1\cup tK_{\ell-1})\\
\ge&(m-\ell+2)\binom{\ell-2}{s-1}+\binom{\ell-2}{s}+\binom{\ell-2}{s-2}-\frac{n_1}{2\ell-1}\binom{2\ell-1}{s}-t\binom{\ell-1}{s}\\
=&\biggl(m-\ell+2\biggr)\binom{\ell-2}{s-1}+\frac{\ell-s-1}{s}\binom{\ell-2}{s-1}+\frac{s-1}{\ell-s}\binom{\ell-2}{s-1}\\
&~~~~~~~~~~~~~~~~~~~~~~~~~~~~-\frac{n_1}{2\ell-1}\binom{2\ell-1}{s}-t\frac{\ell-1}{s}\binom{\ell-2}{s-1}\\
=&\biggl((m-\ell+1)\left(1-\frac{1}{s}\right)+\frac{n_1}{s}+\frac{s-1}{\ell-s}\biggr)\binom{\ell-2}{s-1}-\frac{n_1}{2\ell-1}\binom{2\ell-1}{s}\\
>&\frac{m}{2}\binom{\ell-2}{s-1}-\frac{n_1}{2\ell-1}\binom{2\ell-1}{s}>0.
\end{split}\]
Therefore, the latter part holds.
$\hfill\blacksquare$

\vskip 3mm

\noindent\textbf{Proof of Theorem \ref{Thm2}.} Let $F=P_{\ell_1}\cup\cdots\cup P_{\ell_k}$ be a linear forest with $k\ge2, \ell_i\not=3$ and $G$ an $F$-free graph of order $n$ with $\mathcal{N}_s(G)=ex(n,K_s,F)$ and $s\le \delta_F+1$. Obviously,
\begin{align}\label{3.1}
\mathcal{N}_s(G)\ge max\{\mathcal{N}_s(G_F(n)), ex(n,K_s,P_{\ell_1})\}.
\end{align}
Applying $\left(\delta_F-1\right)$-$disintegration$ to $G$ and let $G^m$ be the $\delta_F$-$core$ of $G$. If $m=0$, then since one vertex is deleted at each step during the process of $(\delta_F-1)$-$disintegration$, that destroys at most $\binom{\delta_F-1}{s-1}$ copies of $K_s$ in $G$, we have
 \[
\begin{split}
\mathcal{N}_s(G)&\le \binom{\delta_F}{s}+\left(n-\delta_F\right) \binom{\delta_F-1}{s-1}
\\&<\binom{\delta_F}{s}+\left(n-\delta_F\right)\binom{\delta_F}{s-1}\le \mathcal{N}_s(G_F(n)),
\end{split}
\]
which contradicts (\ref{3.1}). This implies that $\delta(G^m)\ge \delta_F$ and $m\ge \delta_F+1$.

By the definitions of $\left(\delta_F-1\right)$-$disintegration$, we have
\begin{equation}\label{3.2}
\mathcal{N}_s\left(G^{j+1}\right)-\mathcal{N}_s\left(G^{j}\right)\le \binom{\delta_F-1}{s-1}.
\end{equation}
If  $m\leq j\leq n-1$, then
\begin{equation}\label{3.3}
\mathcal{N}_s(G_F(j+1))-\mathcal{N}_s(G_F(j))\ge \binom{\delta_F}{s-1}.
\end{equation}
Recall that $s\le \delta_F+1$, by (\ref{3.2}) and (\ref{3.3}), we have
$$\mathcal{N}_s\left(G^{j}\right)-\mathcal{N}_s\left(G_F(j)\right)\ge \mathcal{N}_s\left(G^{j+1}\right)-\mathcal{N}_s(G_F(j+1))+1,$$
and so
$$\mathcal{N}_s\left(G^{m}\right)-\mathcal{N}_s\left(G_F(m)\right)\ge \mathcal{N}_s\left(G^{n}\right)-\mathcal{N}_s(G_F(n))+(n-m),$$
which implies that $$\mathcal{N}_s(G^m)\ge n-m+\mathcal{N}_s(G_F(m))\geq n-\delta_F.$$
On the other hand, since $G^m$ has no path of order $|F|$, by Theorem \ref{Luo}, we have
$$ \mathcal{N}_s(G^m)\le \frac{m}{|F|-1}\binom{|F|-1}{s}.$$
Thus, note that  $n\ge 5\binom{|F|-1}{s}^2/\left({(|F|-1)\binom{\delta_F}{s-1}}\right)+\delta_F$, we have

\begin{equation}\label{3.4}
m\ge  \left(n-\delta_F\right)\left(|F|-1\right)\left/\binom{|F|-1}{s}\right.
\ge 5\binom{|F|-1}{s}\left/{\binom{\delta_F}{s-1}}\right.,
\end{equation}
which means the core $G^m$ contains indeed at least $5\binom{|F|-1}{s}/{\binom{\delta_F}{s-1}}$ vertices.

If $G^m$ is disconnected, we assume that the components of $G^m$ are $G_1^m,G_2^m,...,G_c^m$ with $|G_1^m|\geq |G_2^m|\geq \cdots \geq |G_c^m|$.
\vskip 2mm
(i) $F\neq 2P_\ell$ with $\ell$ odd. By (\ref{3.4}), we have $m\ge 3|F|$. If $G^m$ is connected, then by Theorem \ref{Thm1}, $G^m\in \{G_F(m), K_1+tK_{\ell-1},K_2+tK_4\}$. By Claim \ref{c2}, we can see
\[\begin{split}
\mathcal{N}_s(G)&\le \mathcal{N}_s(G^m)+(n-m)\binom{\delta_F-1}{s-1}\le\mathcal{N}_s(G_F(m))+(n-m)\binom{\delta_F-1}{s-1}\\
&\le\mathcal{N}_s(G_F(m))+(n-m)\binom{\delta_F}{s-1}\le\mathcal{N}_s(G_F(n)),
\end{split}\]
with equality holds if and only if $n=m$ and $G=G_F(n)$.

If $G^m$ is disconnected, then since $m\ge 3|F|$, we have $c\ge 4$ by Claim \ref{c1}. If two of the components of $G^m$ are contained in one component of $G$, then we can choose $G_i^m, G_j^m,G_k^m$ such that $G_i^m,G_j^m$ are contained in one component of $G-G_k^m$. Because $\delta(G^m)\geq \delta_F$, we can see that $G-G_k^m$ contains $P_{2\delta_F+2}$ and $G_k^m$ contains $P_{\delta_F+1}$, which implies $G$ contains $F$ by Lemma  \ref{claw2}, a contradiction. Thus, we may assume $G_1,G_2,...,G_s$ are components of $G$ and $G^m_i \subseteq G_i, 1\le i\le c\leq s$. In this case, each $G_i$ contains a $P_{\delta_F+1}$.
If $G$ contains a $P_{\ell_1}$, say $P_{\ell_1}\subseteq G_1$, then since $c\geq 4$, $G$ contains $P_{\ell_1}\cup 3P_{\delta_F+1}$. Assume that $x,y$ are the endpoints of two distinct $P_{\delta_F+1}$. Then $P_{\delta_F+1}\cup \{xy\}\cup P_{\delta_F+1}$ is a path of oder  $2\delta_F+2$, and so $3P_{\delta_F+1}\cup \{xy\}$ contains $F$ by Lemma  \ref{claw2}. If $F$ uses $xy$ to form some $P_{\ell_i}$, then replace this $P_{\ell_i}$ with a $P_{\ell_i}$ from $G_1$, we  find an $F$ in $G$, a contradiction. This is to say that $G$ is $P_{\ell_1}$-free and so
 $\mathcal{N}_s(G)\le ex(n,K_s,P_{\ell_1})$. Therefore, we can obtain that $\mathcal{N}_s(G)\le max\{\mathcal{N}_s(G_F(n)), ex(n,K_s,P_{\ell_1})\}$. Combining this with (\ref{3.1}), we get that  $\mathcal{N}_s(G)= max\{\mathcal{N}_s(G_F(n)), ex(n,K_s,P_{\ell_1})\}$.
\vskip 3mm

(ii)  $F=2P_\ell$ and $\ell$ is odd. By Theorem \ref{Thm2} and Claim \ref{c1}, $G^m$ is $G_F(m)$ or $K_1+(t_1K_{\ell-2}\cup t_2K_{\ell-1})$ or $G^m_1\cup tK_{\ell-1}$. If $s\le \ell-2$,  then because $m\ge 5\binom{|F|-1}{s}/{\binom{\delta_F}{s-1}}$, by Claim \ref{c2},  we have
\[\begin{split}
\mathcal{N}_s(G)&\le \mathcal{N}_s(G^m)+(n-m)\binom{\delta_F-1}{s-1}\\
&\le\mathcal{N}_s(G_F(m))+(n-m)\binom{\delta_F}{s-1}\le\mathcal{N}_s(G_F(n)).
\end{split}\]
At last, if $s=\delta_F+1=\ell-1$, recall that each time we delete a vertex during the process of $(\delta_F-1)$-$disintegration$, no copy of $K_s$ is destroyed. Therefore, $\mathcal{N}_s(G)=\mathcal{N}_s(G^m)$. It is easy to show that $K_1+(tK_{\ell-1}\cup K_r)$ contains the most copies of $K_s$ among all kinds of $\delta_F$-core
when $n$ is large,  where $n-1\equiv r$ (mod $\ell-1$), that is $\mathcal{N}_s(G)=\lfloor\frac{n-1}{\ell-1}\rfloor\ell+1_{r=\ell-2}$.
\vskip 2mm

Combining (i) and (ii), we complete the proof. $\hfill\blacksquare$

\newpage
\noindent{\bf\large Acknowledgements}
\vskip 3mm

This research was supported by NSFC under grant numbers  11871270 and 11931006.

\end{document}